\documentclass[a4paper]{amsart}
\usepackage{amsfonts,amsmath,amssymb,amscd,amstext,amsbsy,latexsym}
\theoremstyle{definition}

\theoremstyle{remark}

\numberwithin{equation}{section}

\begin{document}

\title{A one-line undergraduate proof of Zariski's Lemma and Hilbert's Nullstellensatz}%
\author{Alborz Azarang}%

\maketitle

\centerline{Department of Mathematics, Shahid Chamran University,
Ahvaz, Iran} \centerline{a${}_{-}$azarang@scu.ac.ir}


\vspace{0.5cm}

Our aim in this very short note is to show that the proof of the
following theorem is nothing but the proof of the trivial fact
that $\mathbb Q$ is not a finitely generated $\mathbb
Z$-algebra.\\

{\large\bf Theorem}. {\it Let the domain
$R=K[\alpha_1,\alpha_2,\ldots,\alpha_n]$, where $K$ is field, be
a field. Then every $\alpha_i$ is algebraic over $K$. In
particular, if $K$ is algebraically closed, then $\alpha_i\in K$
for all $i$.}\\

The first part of this theorem is Zariski\rq{}s Lemma and the
second part is usually called Hilbert\rq{}s Nullstellensatz (Weak
Form), which says every maximal ideal $M$ of the polynomial ring
$K[x_1,x_2,\ldots,x_n]$, where $K$ is an algebraically closed
field is of the form $M=(x_1-a_1,x_2-a_2,\ldots,x_n-a_n)$,
$a_i\in K$ for all $i$.

Usually the proof of Zariski\rq{}s Lemma depends on two technical
lemmas due to Artin-Tate and Zariski, see \cite[Proposition 3.2,
and the comment following it]{kunz}. For an elegant proof of
Hilbert\rq{}s Nullstellensatz, based on $G$-ideal theory, see
\cite{kap}. We should also remind the reader that in some
elementary text books on Algebraic Geometry, Noether
normalization theorem is applied for the proof of Zariski\rq{}s
Lemma, see for example \cite[Theorem 1.15]{klus} and \cite{reid}.
Two different proofs of this lemma are also given in \cite[P.
166]{zarsam}. Before presenting our proof, we should emphasize here that, it seems,  the next proof is the simplest and the shortest possible proof (among the existing proofs) of the results in the title, for now.\\

{\bf Proof of the theorem}: Without loss, by induction we may
suppose that $\alpha_1$ is not algebraic over $K$, but
$\alpha_2,\ldots,\alpha_n$ are algebraic over $K(\alpha_1)$ and
get a contradiction (note, R=
$K(\alpha_1)[\alpha_2,\alpha_3,\ldots,\alpha_n]$). Clearly, each
$\alpha_i$ is algebraic over $K[\alpha_1]$. This implies that
there are polynomials $f_2(\alpha_1),\ldots,f_n(\alpha_1)$ in
$K[\alpha_1]$ such that each $\alpha_i$ is integral over the
domain
$A=K[\alpha_1][\frac{1}{f_2(\alpha_1)},\ldots,\frac{1}{f_n(\alpha_1)}]$.
Since $R$ is integral over $A$, we infer that $A$ is a field
(note, it is manifest that if a field $F$ is integral over a
subdomain $D$, then $D$ is a field). Consequently,
$A=K(\alpha_1)$, which is absurd (note, if $D$ is any PID with
infinitely many primes, its field of fractions is never a finitely
generated $D$-algebra). The final part is now
evident.\\

\centerline{{\bf\large Acknowledgment}}

I would like to thank professor O.A.S. Karamzadeh for a useful
discussion on this note and for his encouragement. In particular,
I admit that the following comment of Karamzadeh in his course on
\lq\lq{}Elementary Algebraic Geometry\rq\rq{} which was given 15
years ago, when I was an undergraduate student at Shahid Chamran
University, has always given me a motivation to give a simple
undergraduate proof for the above theorem: To prove Hilbert\rq{}s
Nullstellensatz (Weak Form) in his lecture, he used to prove it
for the maximal ideals in $K[x_1,x_2,\ldots,x_n,\ldots]$, where
$K$ is an uncountable algebraically closed field and also used to
emphasize that the proof in this case is simpler than the proof of
the usual Hilbert\rq{}s Nullstellensatz (Weak Form), see
\cite{karammos}.


\end{document}